\documentclass[12pt,reqno]{amsart}
\usepackage{amssymb}
\usepackage[latin1]{inputenc}
\usepackage{dsfont}
\usepackage{hyperref}
\usepackage{fullpage}
\newtheorem{theorem}{Theorem}[section]
\newtheorem{lemma}[theorem]{Lemma}
\newtheorem{proposition}[theorem]{Proposition}
\newtheorem{corollary}[theorem]{Corollary}
\theoremstyle{definition}
\newtheorem{definition}[theorem]{Definition}

\theoremstyle{remark}
\newtheorem{remark}[theorem]{Remark}
\numberwithin{equation}{section}
\begin{document}
\title{A new kind of augmentation of filtrations}
\author[J. Najnudel]{Joseph Najnudel}
\address{Institut f\"ur Mathematik, Universit\"at Z\"urich, Winterthurerstrasse 190,
8057-Z\"urich, Switzerland}
\email{\href{mailto:joseph.najnudel@math.uzh.ch}{joseph.najnudel@math.uzh.ch}}
\author[A. Nikeghbali]{Ashkan Nikeghbali}
\email{\href{mailto:ashkan.nikeghbali@math.uzh.ch}{ashkan.nikeghbali@math.uzh.ch}}
\date{\today}
\begin{abstract}
Let $(\Omega,\mathcal{F},(\mathcal{F}_t)_{t \geq 0},\mathbb{P})$ be a filtered probability space satisfying the usual assumptions: it is usually not possible to extend to $\mathcal{F}_{\infty}$ 
(the $\sigma$-algebra generated by $(\mathcal{F}_t)_{t \geq 0}$) a coherent family
 of probability measures $(\mathbb{Q}_t)$ indexed by $t \geq 0$, each
of them being defined on $\mathcal{F}_t$. It is known that for instance, on the Wiener space, this extension problem has a positive answer if one takes the filtration generated by the coordinate process,  made right continuous, but can have a negative answer if one takes its usual augmentation. On the other hand, the usual assumptions are crucial in order to obtain the existence of regular versions of paths (typically adapted and continuous or adapted and c\`adl\`ag versions) for most stochastic processes of interest, such as the local time of the standard Brownian motion, stochastic integrals, etc. For instance we shall prove that on the Wiener space, equipped with the right continuous augmentation of the filtration generated by the canonical process, there exists no c\`adl\`ag or continuous and adapted version for the local time at level zero of the canonical process. Hence there is an incompatibility between the problem of extending a coherent family of probability measures to $\mathcal{F}_\infty$ and the classical construction of regular versions of paths, requiring  the usual assumptions (this situation typically occurs in the problem of penalization of the Brownian paths or in mathematical finance), which to the best of our knowledge, has not been noticed so far. 
 In order to fix this problem, we introduce a new property for filtrations, intermediate between the right continuity
and the usual conditions. More precisely, we say that a filtration $(\mathcal{F}_t)_{t \geq 0}$ satisfies the \textit{N-usual} assumptions
if it is right-continuous and if $\mathcal{F}_0$ contains all the sets included in a countable union
of negligible sets $(B_n)_{n \geq 1}$, such that $B_n \in \mathcal{F}_n$ for $n \geq 1$.
There is a natural way to obtain, from a given filtration $(\mathcal{F}_t)_{t \geq 0}$, a new filtration
 which satisfies the N-usual assumptions: we call it the N-augmentation of $(\mathcal{F}_t)_{t \geq 0}$.
We show that most of the important results of the theory of stochastic processes which are generally proved under the usual augmentation, such as the existence of regular version of trajectories or the d\'ebut theorem,   still hold under the N-augmentation; moreover this new augmentation allows the extension of a coherent family of probability measures whenever this is possible with the original filtration. For sake of completeness, we also recall (not so well known) Parthasarathy type conditions on the underlying filtration under which the extension problem for a coherent family of probability measures has a solution. In particular, we shall see that this is always the case on the following two fundamental  spaces: $\mathcal{C}(\mathbb{R}_+,\mathbb{R})$, the space of continuous functions equipped with the filtration generated by the coordinate process and $\mathcal{D}(\mathbb{R}_+,\mathbb{R})$, the space of c\`adl\`ag functions endowed with the filtration generated by the coordinate process.
 \end{abstract}
\maketitle
%\tableofcontents
\section*{Warning}
After we put a preprint version of this paper on the arxiv, it was noticed by Ramon van Handel that the new augmentation we are proposing here and that we call N-augmentation was already introduced by K. Bichteler in 2002 in his book on stochastic integration \cite{B}. K. Bichteler   called this augmentation the natural augmentation. Since the methods we are using are different and since some of our results are not contained in \cite{B}, we have decided to leave our original version unchanged so that the reasons why we were led to naturally introduce this augmentation are still reflected in  the present work. 

\section{Introduction}
\noindent
In stochastic analysis, most of the interesting properties of continuous time random processes cannot be established if 
one does not  assume that their trajectories satisfy some  regularity conditions. For example,
a nonnegative c\`adl\`ag martingale converges almost surely, but if the c\`adl\`ag assumption
is removed, the result becomes false in general. Recall a counter-example: on the filtered
probability space $\big((\mathcal{C}( \mathbb{R}_+, \mathbb{R}), \mathcal{F}, (\mathcal{F}_t)_{t \geq 0}
, \mathbb{W} \big)$, where $\mathcal{F}_t = \sigma \{ X_s, 0 \leq s \leq t\}$, 
$\mathcal{F} = \sigma \{X_s, s \geq 0\}$, $(X_s)_{s \geq 0}$ is the canonical process and $ \mathbb{W}$ the Wiener measure, the 
martingale 
$$\big(M_t := \mathds{1}_{X_t = 1}\big)_{t \geq 0},$$
which is a.s. equal to zero for each fixed $t \geq 0$, does not converge at infinity. 
That is the reason  why one generally  considers a c\`adl\`ag version of  a martingale. 
However there are fundamental examples of stochastic processes for which such a version does not exist. Indeed
let us define on the filtered probability
space $\big((\mathcal{C}( \mathbb{R}_+, \mathbb{R}), \mathcal{F}, (\mathcal{F}_t)_{t \geq 0}
, \mathbb{W} \big)$ described above, the stochastic process $(\mathcal{L}_t)_{t \geq 0}$ as follows:
$$\mathcal{L}_t  = \Phi \left(\underset{m \rightarrow \infty}{\lim \inf} \int_0^t f_m(X_s) ds \right),$$
where $f_m$ denotes the density of a centered Gaussian variable with variance $1/m$ and
$\Phi$ is the function from $\mathbb{R}_+ \cup \{\infty\}$ to $\mathbb{R}_+$ such that 
$\Phi(x) = x$ for $x < \infty$ and $\Phi(\infty) = 0$. The process $(\mathcal{L}_t)_{t \geq 0}$ 
is a version of the local time of the canonical process at level zero, which is defined everywhere
and $(\mathcal{F}_t)_{t \geq 0}$-adapted. It is known that  the process:
$$\big(M_t := |X_t| - \mathcal{L}_t \big)_{t \geq 0}$$
is an $(\mathcal{F}_t)_{t \geq 0}$-martingale. However, $(M_t)_{t \geq 0}$ does not admit 
a c\`adl\`ag version which is adapted. In other words, there exists no c\`adl\`ag, adapted
version $(L_t)_{t \geq 0}$ for the local time at level zero of the canonical process! This 
property can be proved in the following way: let us consider an Ornstein-Uhlenbeck process
 $(U_t)_{t \geq 0}$, starting from zero, and
let us define the process $(V_t)_{t \geq 0}$ by:
$$V_t = (1-t)U_{t/(1-t)}$$
for $t < 1$, and 
$$V_t = 0$$
for $t \geq 1$. This process is a.s. continuous: we denote by $\mathbb{Q}$ its distribution.
One can check the following properties:
\begin{itemize}
\item For all $t \in [0,1)$, the restriction of $\mathbb{Q}$ to $\mathcal{F}_t$ is absolutely continuous
with respect to the corresponding restriction of $\mathbb{W}$. 
\item Under $\mathbb{Q}$, $\mathcal{L}_t\to\infty$  a.s. when $ t\to1,\;t < 1$. 
\end{itemize}
\noindent
By the second property, the set $\{ \mathcal{L}_t \underset{t \rightarrow 1, t<1}{\longrightarrow} \infty\}$
 has probability one under $\mathbb{Q}$. Since it is negligible under $\mathbb{P}$, it is essential 
to suppose that it is not contained in $\mathcal{F}_0$, if we need to have the first property: the filtration
must not to be completed. 
The two properties above imply
\begin{align*}
\mathbb{Q} \left[ L_{1 - 2^{-n}}  \underset{n \rightarrow \infty}{\longrightarrow} \infty \right]
& \geq \mathbb{Q} \left[ \mathcal{L}_{1 - 2^{-n}}  \underset{n \rightarrow \infty}{\longrightarrow} \infty,
\, \forall n \in \mathbb{N}, \, L_{1-2^{-n}} = \mathcal{L}_{1 - 2^{-n}}  \right] \\ 
& \geq 1 - \sum_{n \in \mathbb{N}} \mathbb{Q} \left[L_{1-2^{-n}} \neq \mathcal{L}_{1 - 2^{-n}} \right] = 1. 
\end{align*}
\noindent
The last equality is due to the fact that for all $n \in \mathbb{N}$,
$$ \mathbb{W} \left[L_{1-2^{-n}} \neq \mathcal{L}_{1 - 2^{-n}} \right] = 0,$$
and then
$$\mathbb{Q}  \left[L_{1-2^{-n}} \neq \mathcal{L}_{1 - 2^{-n}} \right] = 0,$$
since $L_{1-2^{-n}}$ and $\mathcal{L}_{1-2^{-n}}$ are $\mathcal{F}_{1-2^{-n}}$-measurable 
and since the restriction of $\mathbb{Q}$ to this $\sigma$-algebra is absolutely continuous with respect to
$\mathbb{W}$. We have thus proved that there exist some paths such that $L_{1 - 2^{-n}}$ tends
to infinity with $n$, which contradicts the fact  that $(L_t)_{t \geq 0}$ is c\`adl\`ag.
From this we also deduce that in general there do not exist c\`adl\`ag versions for martingales. Similarly many other
important results from stochastic analysis cannot be proved
 on the most general filtered probability space, e.g. the existence of the Doob-Meyer decomposition 
for submartinagales and the d\'ebut theorem (see for instance \cite{DM} and \cite{DM2}).
In order to avoid this technical problem, it is generally assumed that the filtered probability 
space on which the processes are constructed satisfies the usual conditions, i.e. the filtration 
is complete and right-continuous. \\

But now, if we wish to perform a change of probability measure (for example,
by using the Girsanov theorem), this assumption reveals to be too restrictive. Let us illustrate this fact by a simple example. Let us consider the filtered probability 
space $\big(\mathcal{C}(\mathbb{R}_+, \mathbb{R}), \widetilde{\mathcal{F}},
 (\widetilde{\mathcal{F}}_t)_{t \geq 0}, \widetilde{\mathbb{W}}\big)$
 obtained, from the Wiener space 
$\big(\mathcal{C}(\mathbb{R}_+, \mathbb{R}), \mathcal{F}, (\mathcal{F}_t)_{t \geq 0}, \mathbb{W}\big)$
described above, by taking its usual augmentation, i.e.:
\begin{itemize}
\item $\widetilde{\mathcal{F}}$ is the $\sigma$-algebra generated by $\mathcal{F}$ and its negligible sets.
\item For all $t \geq 0$, $\widetilde{\mathcal{F}}_t$ is $\sigma$-algebra generated by 
$\mathcal{F}_t$ and the negligible sets of $\mathcal{F}$.
\item $\tilde{\mathbb{W}}$ is the unique possible extension of $\mathbb{W}$ to the completed
$\sigma$-algebra $\widetilde{\mathcal{F}}$.
\end{itemize}
\noindent
Let us also consider the family of probability measures $(\mathbb{Q}_t)_{t \geq 0}$, 
such that $\mathbb{Q}_t$ is defined on $\widetilde{\mathcal{F}}_t$ by
$$\mathbb{Q}_t = e^{X_t - \frac{t}{2} }. \, \widetilde{\mathbb{W}}_{ \, |\widetilde{\mathcal{F}_t} }.$$
This family of probability measures is coherent, i.e. for $0 \leq s \leq t$, the restriction of $\mathbb{Q}_t$ to 
$\widetilde{\mathcal{F}}_s$ is equal to $\mathbb{Q}_s$. However, unlike what one would expect,
 there does not exist a probability measure 
$\mathbb{Q}$ on $\widetilde{\mathcal{F}}$ such that its restriction to 
$\widetilde{\mathcal{F}}_s$ is equal to $\mathbb{Q}_s$ for all $s \geq 0$. Indeed let us assume
that $\mathbb{Q}$ exists. 
The event
$$A := \{ \forall t \geq 0, X_t \geq -1 \} $$
satisfies $\widetilde{\mathbb{W}} [A] = 0$, 
and then $A \in \widetilde{\mathcal{F}}_0$ by completeness, which implies that $\mathbb{Q}[A] = 0$.
On the other hand, under $\mathbb{Q}_t$, for all $t \geq 0$, the process $(X_s)_{0 \leq s \leq t}$ is a Brownian motion
with drift 1, and hence under $\mathbb{Q}$. One
deduces that:
$$\mathbb{Q} [ \forall s \in [0,t], \, X_s \geq -1 ]
= \widetilde{\mathbb{W}} [\forall s \in [0,t], \, X_s \geq -s-1 ] 
\geq \widetilde{\mathbb{W}} [\forall s \geq 0, \, X_s \geq -s-1 ].$$
Consequently by letting $t$ go to infinity one obtains
$$\mathbb{Q} [A] \geq \widetilde{\mathbb{W}} [\forall s \geq 0, \, X_s \geq -s-1 ] > 0,$$
which is a contradiction. Therefore, the usual conditions are not suitable  for the problem of extension of coherent 
probability measures.  In fact one can observe that the argument above does not depend on the completeness of $\widetilde{\mathcal{F}}$, but
only on the fact that $\widetilde{\mathcal{F}}_0$ contains all the sets in $\widetilde{\mathcal{F}}$ of 
probability zero. That is why it still remains available if we consider, with the notation above, the 
space $\big(\mathcal{C}(\mathbb{R}_+, \mathbb{R}), \mathcal{F}, (\mathcal{F}'_t)_{t \geq 0}, \mathbb{W}\big)$,
where for all $t \geq 0$, $\mathcal{F}'_t$ is the $\sigma$-algebra generated by $\mathcal{F}_t$ and 
the sets in $\mathcal{F}$ of probability zero. 

In order to illustrate our point, we now show that if one is not careful with the completion, then one can easily obtain an extension of the Lebesgue measure to all subsets of the real line. Indeed, assume that  
on the measurable space $(\Omega, \widetilde{\mathcal{F}})$ (without any filtration), there 
 exists a probability measure under which the canonical process is a Brownian motion with drift one.
Since the event $\{X_t \underset{t \rightarrow \infty}{\longrightarrow} \infty\}$ 
is in $\mathcal{F}$ and has probability zero under the Wiener measure, all its subsets are contained 
in $\widetilde{\mathcal{F}}$. Then, if under a probability $\mathbb{Q}$, $(X_t)_{t \geq 0}$ is 
a Brownian motion with drift one, one can define a finite measure $\mu$ on $\mathcal{P} ([0,1])$ 
by: 
$$\mu (A) = \mathbb{Q} \left[F(X_1-1) \in A , \, X_t \underset{t \rightarrow \infty}
{\longrightarrow} \infty \right],$$
for all subsets $A$ of $[0,1]$, where $F$ is the distribution function of the standard Gaussian variable. 
Since under $\mathbb{Q}$, $X_t$ tends almost surely to infinity when $t\to\infty$, one checks that 
$\mu$ is an extension of the Lebesgue measure, which is defined for all the subsets of $[0,1]$. Now, the 
existence of such an extension is incompatible with the continuum hypothesis if 
one assumes the usual axioms of set theory (see \cite{cies} for details and references on this problem).

The above discussion (which corresponds to the simplest non-trivial Girsanov transformation) outlines that
by doing the usual augmentation, we add, in a certain sense, too much information in the filtration which is 
considered: in particular the asymptotic properties of the processes are put directly into $\mathcal{F}_0$. It is precisely the fact that asymptotic events are already in  $\mathcal{F}_0$ that it is not possible in general to construct a probability measure $\mathbb{Q}$ which is singular with repect to $\mathbb{W}$ but locally absolutely continuous with respect to it.
Let us shortly sum up the two problems we have encountered emphasized:
\begin{itemize}
\item if we do not complete the filtrations, most of the properties of the trajectories of the stochastic 
processes are lost in general;
\item if we add all the negligible sets in $\mathcal{F}_0$, it is not possible in general to extend to
 $\mathcal{F}_\infty=\sigma(\bigcup_{t\geq0} \mathcal{F}_t)$
 a coherent family of probability measures defined on $\mathcal{F}_t$. In particular, it is not possible to construct, by
 Girsanov transformation, a probability which is singular with respect to the initial probability measure. 
\end{itemize}
\noindent
Now a natural question is: what can one do  to avoid simultaneously these two issues?
This case arises, for instance, if one  needs to perform Giranov transformations with a density involving 
local times. It also occurs when, working with an infinite time horizon,  one defines locally on each $\mathcal{F}_t$ a new probability measure with the help of a martingale which is not uniformly integrable (e.g. $\exp(\gamma X_t-\gamma^2 t/2)$ on the Wiener space) and then one manipulates stochastic integrals, local times, etc. This situation is very often encountered in some problems related to the penalization of the Wiener paths or sometimes in financial modeling.  To the best of our knowledge this important issue has never been dealt with before, hence making the results obtained when facing these two problems not totally rigorous.\\

 The goal of this paper is to show that by performing a new kind of augmentation of filtrations, 
intermediate between the right-continuous version and the usual augmentation,  it is  possible to address the two technical issues addressed above.  This augmentation, called \textit{N-augmentation} in this 
paper, is obtained as follows: after making the filtration right-continuous, instead of putting
all the negligible sets of $\mathcal{F}$ in $\mathcal{F}_0$ (which adds too much information), we only put the 
sets which are contained in a countable union $(B_n)_{n \geq 0}$ of sets of probability zero, such that
$B_n \in \mathcal{F}_n$ for all $n \geq 0$. Note that this way of completing filtrations has a
concrete interpretation, which can be interesting in modeling problems:
the events of probability zero can be anticipated if it concerns the future only up to a finite time, even
if this time is unbounded. 

More precisely the paper is organized as follows. 
\begin{itemize}
\item In Section \ref{definitions}, we construct in details 
the N-augmentation of a filtration and compare it with the usual augmentation;
\item In Section \ref{theorems}, we show that most of the classical theorems of stochastic analysis which are 
proved under the usual conditions remain true under N-usual conditions;
\item In Section \ref{extensions} we show that the N-augmentation preserves the main properties of 
the stochastic processes (martingale property for example), and allows us to deal
 with measures which are singular with respect to the initial probability measure.  
In particular, we review sufficient conditions of Parthasrarthy type under 
which a coherent family of probability measures can be extended, and we show that in the most interesting 
cases, this extension remains possible after taking the N-augmentation. 
\end{itemize}

\section*{Acknowledgements}
We would like to thank Freddy Delbaen for all his enlightening comments during the preparation of this paper, and Hans F\"ollmer for very helpful discussions. We also wish to thank Ramon van Handel for bringing to our attention the book by K. Bichteler where the new augmentation we suggest was first proposed.

\section{The N-augmentation of filtrations} \label{definitions}
\noindent
Let $(\Omega, \mathcal{F}, (\mathcal{F}_t)_{t \geq 0}, \mathbb{P})$ be a filtered probability space.
We first introduce a few definitions in order to  rigorously define  the N-usual conditions.  The first notion we introduce is the notion of N-negligible sets, which corresponds to the 
sets we want to put into $\mathcal{F}_0$ in the N-augmentation. 
\begin{definition} \label{negligible}
A subset $A$ of $\Omega$ is N-negligible with respect to the space 
$(\Omega,\mathcal{F}, (\mathcal{F}_t)_{t \geq 0}, \mathbb{P})$, iff there exists a sequence $(B_n)_{n \geq 0}$
of subsets of $\Omega$, such that for all $n \geq 0$, $B_n \in \mathcal{F}_n$, $\mathbb{P} [B_n] = 0$, and
$$A \subset \bigcup_{n \geq 0} B_n.$$  
\end{definition}
\begin{remark}
The integers do not play a crucial r\^ole in Definition \ref{negligible}. If 
$(t_n)_{n \geq 0}$ is an unbounded sequence in $\mathbb{R}_+$, one
can replace the condition $B_n \in \mathcal{F}_n$ by the condition $B_n \in \mathcal{F}_{t_n}$.  
\end{remark}
\noindent
Let us now define a notion which is the analog of completeness for N-negligible sets. It 
is the main ingredient in the definition of what we shall call the N-usual conditions:  
\begin{definition} \label{complete}
A filtered probability space $(\Omega,\mathcal{F}, (\mathcal{F}_t)_{t \geq 0}, \mathbb{P})$,
  is N-complete iff all the N-negligible sets of this space are contained in $\mathcal{F}_0$. It
satisfies the N-usual conditions iff it is N-complete and the filtration 
$(\mathcal{F}_t)_{t \geq 0}$ is right-continuous.
\end{definition}
\noindent
It is natural to ask if from a given filtered probability 
space, one can define in a canonical way a space which satisfies the N-usual conditions and 
which is as "close" as possible to the initial space. The answer to this question is positive in the following
sense: 
\begin{proposition}
Let $(\Omega, \mathcal{F}, (\mathcal{F}_t)_{t \geq 0}, \mathbb{P})$ be a filtered probability space, and 
$\mathcal{N}$ the family of its N-negligible sets. 
Let $\widetilde{\mathcal{F}}$ be the $\sigma$-algebra generated by $\mathcal{N}$ and 
$\mathcal{F}$, and for all $t \geq 0$, 
 $\widetilde{\mathcal{F}}_t$  the $\sigma$-algebra generated by $\mathcal{N}$ and
$\mathcal{F}_{t+}$,
where $$ \mathcal{F}_{t+} := \bigcap_{u > t} \mathcal{F}_t.$$ Then there exists a unique probability measure 
$\widetilde{\mathbb{P}}$ on $(\Omega, \widetilde{\mathcal{F}})$ which coincides with 
$\mathbb{P}$ on $\mathcal{F}$, and the space $(\Omega, \widetilde{\mathcal{F}}, 
 (\widetilde{\mathcal{F}}_t)_{t \geq 0},
\widetilde{\mathbb{P}})$ satisfies the N-usual conditions. Moreover, if 
 $(\Omega, \mathcal{F}', (\mathcal{F}'_t)_{t \geq 0}, \mathbb{P}')$ is a filtered probability space
satisfying the N-usual conditions, such that $\mathcal{F}'$ contains $\mathcal{F}$,  $\mathcal{F}'_t$ contains
$\mathcal{F}_t$ for all $t \geq 0$, and if $\mathbb{P}'$ is an extension of $\mathbb{P}$, then 
$\mathcal{F}'$ contains $\widetilde{\mathcal{F}}$, $\mathcal{F}'$ contains $\widetilde{\mathcal{F}}_t$, 
for all $t \geq 0$ and $\mathbb{P}'$ is an extension of $\widetilde{\mathbb{P}}$.
In other words,  $(\Omega, \widetilde{\mathcal{F}}, 
 (\widetilde{\mathcal{F}}_t)_{t \geq 0},
\widetilde{\mathbb{P}})$ is the smallest extension of $(\Omega, \mathcal{F}, 
(\mathcal{F}_t)_{t \geq 0}, \mathbb{P})$ which satisfies the N-usual conditions:
we call it the N-augmentation of $(\Omega, \mathcal{F}, 
(\mathcal{F}_t)_{t \geq 0}, \mathbb{P})$ 
\end{proposition}
\begin{proof}
Let us first denote by $\mathcal{E}$ the family of subsets $A$ of $\Omega$ such that 
there exists $A' \in \mathcal{F}$, satisfying:
\begin{equation}
(A \backslash A') \cup (A' \backslash A) \, \in \, \mathcal{N}. \label{symdif}
\end{equation}
If $A \in \mathcal{E}$, if $A'$ satisfies \eqref{symdif}, and if we denote by $B$ the 
complement of $A$, $B'$ the complement of $A'$, then $B' \in \mathcal{F}$ and:
$$(B \backslash B') \cup (B' \backslash B) \, = (A' \backslash A) \cup (A \backslash A') \, \in
\, \mathcal{N},$$
which implies that $B \in \mathcal{E}$. Moreover, if $A_n \in \mathcal{E}$ for all $n \geq 1$, 
then there exists, for all $n \geq 1$, $A'_n \in \mathcal{F}$ such that:
$$ (A_n \backslash A'_n) \cup (A'_n \backslash A_n) \, \in \, \mathcal{N}.$$
One has:
$$\left[\left(\bigcup_{n \geq 1} A_n \right) \, \backslash \, \left( \bigcup_{n \geq 1} A'_n
\right)\right] \,  \cup  \, \left[\left(\bigcup_{n \geq 1} A'_n \right) \, \backslash \, \left( \bigcup_{n \geq 1} 
A_n \right)\right] \, \subset \bigcup_{n \geq 1} \left[(A_n \backslash A'_n) \cup (A'_n \backslash A_n)\right]$$
and then,
$$\left[\left(\bigcup_{n \geq 1} A_n \right) \, \backslash \, \left( \bigcup_{n \geq 1} A'_n
\right)\right] \,  \cup  \, \left[\left(\bigcup_{n \geq 1} A'_n \right) \, \backslash \, \left( \bigcup_{n \geq 1} 
A_n \right)\right] \, \in \, \mathcal{N},$$
since $\mathcal{N}$ is stable by countable union.
Therefore, 
$$ \bigcup_{n \geq 1} A_n \in \mathcal{E},$$
and $\mathcal{E}$ is a $\sigma$-algebra. Since it obviously contains $\mathcal{F}$ and $\mathcal{N}$, 
it contains $\widetilde{\mathcal{F}}$. On the other hand, if $A \in \mathcal{E}$, and 
$A'$ satisfies \eqref{symdif}, then:
$$A = \left[A' \cup (A \backslash A') \right] \, \backslash \, (A' \backslash A) \in \widetilde{\mathcal{F}},$$
since $A' \in \mathcal{F}$, $A \backslash A' \in \mathcal{N}$ and $A' \backslash A \in \mathcal{N}$.
In other words, we have proved:
$$\mathcal{E} = \widetilde{\mathcal{F}}.$$
Similarly, for all $t \geq 0$,  a set $A$ is in $\widetilde{\mathcal{F}}_t$ iff there exists 
$A' \in \mathcal{F}_{t+}$ satisfying \eqref{symdif}. 
Now, let $\widetilde{\mathbb{P}}$ be a probability on $(\Omega, \widetilde{\mathcal{F}})$ extending
$\mathbb{P}$. If $A \in \widetilde{\mathcal{F}}$, and $A' \in \mathcal{F}$ satisfies 
\eqref{symdif}, then
$$\widetilde{\mathbb{P}} [A \backslash A'] = \widetilde{\mathbb{P}} [A' \backslash A] = 0,$$
since these two sets, N-negligible, are included in a set in $\mathcal{F}$ of probability zero.
This implies:
$$\widetilde{\mathbb{P}} [A] = \mathbb{P}[A'] + \widetilde{\mathbb{P}} [A \backslash A'] -
\widetilde{\mathbb{P}} [A' \backslash A] = \mathbb{P} [A']$$ 
and uniqueness of $\widetilde{\mathbb{P}}$ if it exists. To prove existence, let us first observe that 
for three sets $A \in \widetilde{\mathcal{F}}$, $A', A'' \in \mathcal{F}$ satisfying \eqref{symdif}
and the similar equation with $A'$ replaced by $A''$, 
$$\mathbb{P} [A'] = \mathbb{P}[A''].$$
This is a consequence of the fact that $A' \backslash A''$ and $A'' \backslash A'$ are 
in the intersection of $\mathcal{N}$ and $\mathcal{F}$, and then, have probability zero. 
In other words, if $A \in \widetilde{\mathcal{F}}$, one can define:
$$\widetilde{\mathbb{P}} [A] := \mathbb{P} [A'],$$
since $\mathbb{P}[A']$ does not depend of the choice of $A'$ satisfying \eqref{symdif}. 
Since for $A \in \mathcal{F}$, one can take $A'=A$, $\widetilde{\mathbb{P}}$ is an extension of $\mathbb{P}$.
Now let $(A_n)_{n \geq 1}$ be a sequence of disjoint sets in $\widetilde{\mathcal{F}}$, and 
$(A'_n)_{n \geq 1}$ a sequence of sets in $\mathcal{F}$ satisfying the analog of \eqref{symdif}. 
One has
$$\left[\left(\bigcup_{n \geq 1} A_n \right) \, \backslash \, \left( \bigcup_{n \geq 1} A'_n
\right)\right] \,  \cup  \, \left[\left(\bigcup_{n \geq 1} A'_n \right) \, \backslash \, \left( \bigcup_{n \geq 1} 
A_n \right)\right] \, \in \, \mathcal{N},$$
and then
$$\widetilde{\mathbb{P}}
 \left[ \bigcup_{n \geq 1} A_n  \right] = \mathbb{P} \left[
\bigcup_{n \geq 1} A'_n  \right].$$
Now for $1 \leq m \leq n$,
$$\mathbb{P} [A'_m \cap A'_n] = 0,$$
since $A_m$ and $A_n$ are disjoint, and then
$$\widetilde{\mathbb{P}}
 \left[ \bigcup_{n \geq 1} A_n  \right] = \sum_{n \geq 1} \mathbb{P} [A'_n] = \sum_{n \geq 1} 
\widetilde{\mathbb{P}} [A_n],$$
which implies that $\widetilde{\mathbb{P}}$ is a probability measure. Let us now prove that
$(\Omega, \widetilde{\mathcal{F}},  (\widetilde{\mathcal{F}}_t)_{t \geq 0},
\widetilde{\mathbb{P}})$ satisfies the N-usual conditions. The N-completeness can be checked
 in the following way. If
$A$ is an N-negligible set with respect to $(\Omega, \widetilde{\mathcal{F}}, 
 (\widetilde{\mathcal{F}}_t)_{t \geq 0}, \widetilde{\mathbb{P}})$ , there exists $(B_n)_{n \geq 1}$ such that 
$B_n \in \widetilde{\mathcal{F}}_n$ and $\widetilde{\mathbb{P}} [B_n] = 0$ for all $n \geq 1$, and:
$$A \subset \bigcup_{n \geq 1} B_n.$$
Since $B_n \in \widetilde{\mathcal{F}}_n$ there exists $B'_n \in \mathcal{F}_{n+}$ satisfying the analog 
of \eqref{symdif}, which implies
 $$\mathbb{P} [B'_n] = \widetilde{\mathbb{P}} [B_n] = 0,$$ 
and then $B'_n \in \mathcal{N}$. Since $B_n \backslash B'_n \in \mathcal{N}$, $B_n \cup B'_n$ and 
then $B_n$ are in $\mathcal{N}$. Finally the union of $B_n$ for $n \geq 1$ is also in $\mathcal{N}$, 
which implies $A \in \mathcal{N}$ and then $A \in \widetilde{\mathcal{F}}_0$. Let us now prove
the right-continuity of $(\widetilde{\mathcal{F}}_t)_{t \geq 0}$.
Let $t \geq 0$ and let $A$ be in the intersection of $\widetilde{\mathcal{F}}_s$ for $s > t$. 
For all integers $n \geq 1$ there exists $A'_n \in \mathcal{F}_{(t + 1/n)+}$ such that
its symmetric difference with $A$ is in $\mathcal{N}$. One deduces that the symmetric difference 
between $A$ and
$$ A':= \bigcap_{m \geq 1} \bigcup_{n \geq m} A'_n$$
is in $\mathcal{N}$. Now since for all integers $m_0 \geq 1$
$$A' = \bigcap_{m \geq m_0} \bigcup_{n \geq m} A'_n,$$
$A' \in \mathcal{F}_{(t + 1/m_0)+}$ for all $m_0 \geq 1$, which implies that $A' \in \mathcal{F}_{t+}$, 
and then $A \in \widetilde{\mathcal{F}}_t$. 
We have thus proved that $(\Omega, \widetilde{\mathcal{F}},  (\widetilde{\mathcal{F}}_t)_{t \geq 0},
\widetilde{\mathbb{P}})$ satisfies the N-usual conditions; it remains to show that
it is the smallest extension of $(\Omega, \mathcal{F}, (\mathcal{F}_t)_{t \geq 0}, \mathbb{P})$
which enjoys this property. Let $(\Omega, \mathcal{F}', (\mathcal{F}'_t)_{t \geq 0}, 
\mathbb{P}')$ be such an extension. By N-completeness $\mathcal{F}'_0$ contains
all the N-negligible sets of $(\Omega, \mathcal{F}', (\mathcal{F}'_t)_{t \geq 0}, 
\mathbb{P}')$ and a fortiori the N-negligible sets of 
$(\Omega, \mathcal{F}, (\mathcal{F}_t)_{t \geq 0}, \mathbb{P})$. Moreover, for all $t \geq 0$, 
$\mathcal{F}'_t = \mathcal{F}'_{t+}$ by right-continuity, which implies that $\mathcal{F}'_t$ contains
$\mathcal{F}_{t+}$. Since it contains $\mathcal{N}$, it also contains $\widetilde{\mathcal{F}}_t$, and 
similarly, $\mathcal{F}'$ contains $\widetilde{\mathcal{F}}$. Now, since $\mathbb{P}'$ is an 
extension of $\mathbb{P}$, its restriction to $\widetilde{\mathcal{F}}$ is also an extension of $\mathbb{P}$ and 
by uniqueness it is necessarily equal to $\widetilde{\mathbb{P}}$. Hence, $\mathbb{P}'$ is 
an extension of $\widetilde{\mathbb{P}}$.  
\end{proof}
\noindent
Once the N-usual conditions are defined, it is natural to compare them with the usual conditions. 
One has the following result:
\begin{proposition}
Let $(\Omega, \mathcal{F}, (\mathcal{F}_t)_{t \geq 0}, \mathbb{P})$ be a filtered probability
space which satisfies the N-usual conditions. Then for all $t \geq 0$, the space 
$(\Omega, \mathcal{F}_t, (\mathcal{F}_s)_{0 \leq 
s \leq t}, \mathbb{P})$ satisfies the usual conditions. 
\end{proposition}
\begin{proof}
The right-continuity of $(\mathcal{F}_s)_{0 \leq s \leq t}$ is obvious, let us prove the completeness. 
If $A$ is a negligible set of $(\Omega, \mathcal{F}_t, \mathbb{P})$, there exists 
$B \in \mathcal{F}_t$, such that 
$A \subset B$ and $\mathbb{P}[B]=0$. One deduces immediately that $A$ is N-negligible with respect 
to $(\Omega, \mathcal{F}, (\mathcal{F}_t)_{t \geq 0}, \mathbb{P})$, and by N-completeness
of this filtered probability space, $A \in \mathcal{F}_0$. 
\end{proof}
\noindent
This relation between usual conditions and N-usual conditions is the main ingredient of 
the results in Section \ref{theorems}  where we prove that one can replace the usual conditions
by the N-usual conditions in most of the classical results in stochastic calculus.

\section{Classical theorems under the N-usual conditions} \label{theorems}
In the introduction we observed that it is very useful to have c\`adl\`ag versions of
martingales. These versions always exist under the N-usual conditions:
\begin{proposition}
Let $(X_t)_{t \geq 0}$ be a submartingale or a supermartingale, with respect to a filtered
probability space satisfying the N-usual conditions. If $\mathbb{E}[X_t]$ is right-continuous 
with respect to $t$ (in particular if $(X_t)_{t \geq 0}$ is 
a martingale), then $(X_t)_{t \geq 0}$ admits a c\`adl\`ag modification, which is unique 
up to indistinguishability. 
\end{proposition}
\begin{proof}
Let us assume that $(X_t)_{t \geq 0}$ is defined on the filtered probability space
$(\Omega, \mathcal{F}, (\mathcal{F}_t)_{t \geq 0}, \mathbb{P})$. 
For all $t \geq 0$, the process $(X_s)_{0 \leq s \leq t}$ is a submartingale or a supermartingale with 
respect to $(\Omega, \mathcal{F}_t, (\mathcal{F}_s)_{0 \leq s \leq t}, \mathbb{P})$ which satisfies the
usual conditions since $(\Omega, \mathcal{F}, (\mathcal{F}_t)_{t \geq 0}, \mathbb{P})$ satisfies the
N-usual conditions. By right-continuity of $\mathbb{E}[X_t]$ with respect to $t$, $(X_s)_{0 \leq s \leq t}$
admits a c\`adl\`ag modification $(\tilde{X}^{(t)}_s)_{0 \leq s \leq t}$, which is unique up to 
indistinguishability. This uniqueness implies that for $0 \leq t \leq u$, one has almost surely
$$\tilde{X}^{(t)}_s = \tilde{X}^{(u)}_s$$
for all $s \leq t$. Let us denote by $N$ the set of $\omega \in \Omega$ such that
there exists integers $n \geq m \geq 1$, and $s \in [0, m]$, such that 
$$\tilde{X}^{(n)}_s (\omega) \neq \tilde{X}^{(m)}_s (\omega).$$
It is easy to check that $N$ is N-negligible, and then $N \in \mathcal{F}_0$ with 
$\mathbb{P} [N]= 0$.
One can now define a process $(\tilde{X}_s)_{s \geq 1}$ by 
$$\tilde{X}_s = \tilde{X}^{(n)}_s$$
for $n \geq s \geq 0$, on the complement of $N$, and by $\tilde{X}_s = 0$ for all
 $s \geq 0$, on $N$. This process is c\`adl\`ag, adapted (recall that $N \in \mathcal{F}_0$) and 
is a modification of $(X_s)_{s \geq 0}$. The uniqueness of this c\`adl\`ag modification is immediate.
\end{proof}
\noindent
There are many situations  where one can construct c\`adl\`ag, or even continuous, versions
of stochastic processes, but  in general one cannot make sure that these versions are adapted to 
the natural filtration of the initial process. This problem can be solved by introducing the N-augmentation. 
Before stating the corresponding proposition let us prove the following useful lemma:
\begin{lemma} \label{cadlag}
Let $A > 0$, and let $f$ be a real valued function defined on a dense subset $D$ of $[0,A]$, containing $A$. Then 
$f$ can be extended to a c\`adl\`ag function from $[0,A]$ to $\mathbb{R}$ iff the  following two
conditions hold:
\begin{itemize}
\item For all $x \in D \cap [0,A)$, $f(y)$ tends to $f(x)$ for $y \in D$ going to $x$ from above.
\item For all $\epsilon > 0$, there exists $\delta > 0$ such that for all $x, y, z \in D$, 
$x \leq y \leq z \leq x+ \delta$ implies $|f(y)-f(x)| \wedge |f(y)-f(z)|  \leq \epsilon$.
\end{itemize}
\noindent
If these conditions hold, the c\`adl\`ag extension of $f$ is unique.
\end{lemma}
\begin{proof}
If the first condition given above is not satisfied, it is obvious that $f$ cannot be extended to a 
right-continuous function. Let us now suppose that the second condition does not hold. There exist
 $\epsilon >0$ and three sequences $(x_n)_{n \geq 0}$, $(y_n)_{n \geq 0}$, $(z_n)_{n \geq 0}$ 
such that $x_n\leq y_n \leq z_n$
for all $n \geq 0$, $z_n - x_n$ tends to zero when $n$ goes to infinity, and 
$|f(x_n) - f(y_n)|  \wedge |f(y_n)-f(z_n)| > \epsilon$ for all $n \geq 0$. 
By taking subsequences, one can assume that $(y_n)_{n \geq 0}$ 
is monotone, and then converges to a limit $a \in [0,A]$.
One has three possible cases:
\begin{itemize}
\item If $(y_n)_{n \geq 0}$ is non-decreasing, and $y_n < a$ for all $n \geq 0$, then
$x_n < a$. Since $0 \leq y_n-x_n \leq z_n - x_n$ tends to zero when $n$ goes to infinity,
$x_n$ tends to $a$, as $y_n$. We have proved that the two sequences $(x_n)_{n \geq 0}$ and $(y_n)_{n \geq 0}$ 
tend to $a$, strictly from below, but since $|f(x_n)- f(y_n)| > \epsilon$ for all $n$, $f$
 cannot have a left limit at $a$. 
\item If $(y_n)_{n \geq 0}$ is non-decreasing, and $y_{n_0} = a$ for some integer $n_0$, then
$y_n = a$ for all $n \geq n_0$, which implies $z_n \geq a$. Since $0 \leq z_n - y_n \leq z_n - x_n$ 
tends to zero when $n$ goes to infinity, $z_n$ tends to $a$. Hence $(z_n)_{n \geq 0}$ tend to $a$ from above, and since $|f(z_n) - f(y_n)| > \epsilon$ for all $n$,
$|f(z_n)- f(a)| > \epsilon$ for all $n \geq n_0$, and $f(a)$ cannot be the right limit of $f$ at $a$.
\item If $(y_n)_{n \geq 0}$ is non-increasing, $z_n \geq y_n \geq a$ for all $n \geq 0$, and
$z_n$ tends to $a$ when $n$ goes to infinity. Hence, $(y_n)_{n \geq 0}$ and $(z_n)_{n \geq 0}$
tend to $a$ from above. Since $|f(z_n) -f(y_n)| > \epsilon$ for all $n$, $f$ has no right limit
at $a$ if $a \notin D$, and $f(a)$ cannot be the right limit of $f$ at $a$ if $a \in D$. 
\end{itemize}
\noindent
In any case $f$ cannot be extended to a c\`adl\`ag function. 
On the other hand, let us suppose that $f$ has not a right limit at $a \in [0,A)$. 
Then, one can find a sequence $(t_n)_{n \geq 0}$, strictly decreasing to $a$, such that
one of these three properties holds:
\begin{itemize}
\item $f(t_{n+1}) \geq f(t_n) +1$ for all $n \geq 0$.
\item  $f(t_{n+1}) \leq f(t_n) -1$ for all $n \geq 0$.
\item There exist $u, v \in \mathbb{R}$ such that for all $n \geq 0$, 
$f(t_{2n}) < u < v < f(t_{2n+1})$.
\end{itemize}
\noindent
In any case, $t_{2n+2} \leq t_{2n+1} \leq t_{2n}$ for all $n \geq 0$, $t_{2n}-t_{2n+2}$ tends to zero 
when $n$ goes to infinity, but $|f(t_{2n})- f(t_{2n+1})| \wedge |f(t_{2n+1})- f(t_{2n+2})| \geq \epsilon$
for some $\epsilon > 0$, independent of $n$. Hence the second condition given in Lemma \ref{cadlag} is not 
satisfied. One has a similar results if $f$ has no left limit at $a \in (0,A]$. 
Now let us assume that $f$ satisfies the two conditions given in Lemma \ref{cadlag}. Necessarily, 
$f$ admits left and right limits everywhere. Let $g$ be the function from $[0,A]$ to $\mathbb{R}$
such that:
\begin{itemize}
\item For $t < A$, $g(t)$ is the right limit of $f$ at $t$;
\item For $t = A$, $g(t) = f(A)$.
\end{itemize}
\noindent
By assumption, for $t \in D \cap [0,A)$, the right limit of $f$ at $t$ is $f(t)$, and then 
$g$ coincides with $f$ on $D$. It remains to prove that $g$ is c\`adl\`ag. 
Let $t \in [0,A)$. For all $\epsilon > 0$, there exists $u > t$ such that 
$|f(v)- g(t)| \leq \epsilon$ for all $v \in D \cap (t,u)$. Now for all $w \in [t,u)$, 
$g(w)$ is the limit of $f(v)$ for $v \in D \cap (t,u)$ strictly decreasing to $w$, 
which implies that $|g(w)-g(t)| \leq \epsilon$. Hence, $g$ is right-continuous. 
Now, let $t \in (0,A]$, and let $b$ be the left limit of $f$ at $t$. For all $\epsilon > 0$, there 
exists $u< t$ such that $|f(v) - b| \leq \epsilon$ for all $v \in D \cap (u,t)$. For 
$w \in (u,t)$, $g(w)$ is the limit of $f(v)$ for $v \in D \cap (u,t)$ strictly decreasing 
to $w$, and then $|g(w) - b| \leq \epsilon$. Therefore, $g$ admits left limits.
The uniqueness of the c\`adl\`ag extension of $f$ is due to the fact that two c\`adl\`ag functions
from $[0,A]$ to $\mathbb{R}$ which coincide on a dense subset of $[0,A]$ containing $A$ are 
necessarily equal.
\end{proof}
\noindent
We are now able to prove the following result:
\begin{proposition} \label{cadlagversion}
Let $(\Omega, \mathcal{F}, (\mathcal{F}_t)_{t \geq 0}, \mathbb{P})$ be a filtered probability space
satisfying the N-usual conditions, and let $(X_t)_{t \geq 0}$ be an adapted process defined on this space. 
We assume that there exists a c\`adl\`ag version $(Y_t)_{t \geq 0}$
 of $(X_t)_{t \geq 0}$. Then there exists a c\`adl\`ag and adapted version of 
$(X_t)_{t \geq 0}$, which is necessarily indistinguishable from $(Y_t)_{t \geq 0}$. 
\end{proposition}
\begin{proof}
Let $D$ be a countable and dense subset of $\mathbb{R}_+$, containg $\mathbb{N}$,
 and for all integers $n \geq 1$, let $N_n$ be  
the set of $\omega \in \Omega$ such that the function $f$ from $D \cap [0,n]$ to $\mathbb{R}$, 
defined by $f(t) = X_t(\omega)$, does not admit a unique c\`adl\`ag extension to $[0,n]$.
By Lemma \ref{cadlag}, $N_n \in \mathcal{F}_n$ for all $n \geq 0$, 
since $(X_t)_{t \geq 0}$ is adapted. 
Now, $(Y_t)_{t \geq 0}$ is c\`adl\`ag, which implies that for all $n$ 
$$N_n \subset \{\exists t \in D, X_t \neq Y_t\}.$$
Since $(Y_t)_{t \geq 0}$ is a version of $(X_t)_{t \geq 0}$, one deduces that 
$\mathbb{P} [N_n] = 0$ for all $n$. Hence
 $$N := \bigcup_{n \geq 0} N_n$$
is N-negligible, which implies that $N \in \mathcal{F}_0$ and $\mathbb{P}[N]=0$.
Now, let $\omega \notin N$, and for $n \geq 1$, let $g_{\omega,n}$ be the unique c\`adl\`ag extension to $[0,n]$
of the function $f_{\omega,n}$ from $D \cap [0,n]$ to $\mathbb{R}$, 
defined by $f_{\omega,n}(t) = X_t(\omega)$. By uniqueness, $g_{\omega,m}$ and $g_{\omega,n}$ coincide 
on $[0,m]$ for $m \leq n$. Hence there exists a c\`adl\`ag function $g_{\omega}$ from 
$\mathbb{R}_+$ to $\mathbb{R}$ such that $g_{\omega}(t) = X_t(\omega)$ for all $t \in D$. 
Now, let $(\tilde{X}_t)_{t \geq 0}$ be the process defined by
$$\tilde{X}_t( \omega) = g_{\omega} (t) \, \mathds{1}_{\omega \notin N},$$
which is c\`adl\`ag (on $N$ it is identically zero). For all $t \geq 0$, and 
for any sequence $(t_n)_{n \geq 0}$ of elements of $D$, tending to $t$ from above
$$\tilde{X}_t(\omega) = \mathds{1}_{\omega \notin N} \, \underset{n \rightarrow \infty}{\lim} 
g_{\omega}(t_n) = \mathds{1}_{\omega \notin N} \, \underset{n \rightarrow \infty}{\lim} 
X_{t_n} (\omega).$$
Since $(X_t)_{t \geq 0}$ is adapted, $(\mathcal{F}_t)_{t \geq 0}$ is right-continuous and 
$N \in \mathcal{F}_0$, $(\tilde{X}_t)_{t \geq 0}$ is adapted. 
Moreover for all $t \geq 0$, $\tilde{X}_t$ is almost surely the right limit of $Y$ at $t$, restricted to $D$, 
since $Y$ is a version of $X$ and $N$ is negligible. Since $Y$ is c\`adl\`ag, $\tilde{X}_t = Y_t$ a.s., 
and $\tilde{X}_t = X_t$ almost surely. Consequently $\tilde{X}$ is a c\`al\`ag and adapted version of $X$. 
Now two c\`adl\`ag versions of $X$ are necessarily indistinguishable, since there almost surely coincide
at all $t \in D$.
\end{proof}
One has a similar result for continuous processes:
\begin{proposition} 
Let $(\Omega, \mathcal{F}, (\mathcal{F}_t)_{t \geq 0}, \mathbb{P})$ be a filtered probability space
satisfying the N-usual conditions, and let $(X_t)_{t \geq 0}$ be an adapted process defined on this space. 
We assume that there exists a continuous version $(Y_t)_{t \geq 0}$
 of $(X_t)_{t \geq 0}$. Then there exists a continuous and adapted version of 
$(X_t)_{t \geq 0}$, which is necessarily indistinguishable from $(Y_t)_{t \geq 0}$. 
\end{proposition}
\begin{proof}
It is similar to the proof of Proposition \ref{cadlagversion}, and more simple, so we go quickly. 
If $D$ is a countable and dense subset of $\mathbb{R}_+$, the restriction of $X$ to $D$ is,
except on an N-negligible set $N$, locally,
uniformly continuous, since $X$ has a continuous version. Therefore one can define a continuous, adapted version
of $X$ by taking its limit (after restriction to $D$), on the complement of $N$, and zero on $N$. 
\end{proof}
\noindent
Similarly one has the following version of the Kolmogorov criteria:
\begin{proposition}
Let $(\Omega, \mathcal{F}, (\mathcal{F}_t)_{t \geq 0}, \mathbb{P})$ be a filtered probability space 
satisfying the N-usual conditions, and let $(X_t)_{t \geq 0}$ be an adapted process defined on this space. 
We suppose that there exist $\alpha, \beta > 0$, and for all $A>0$, $C(A) > 0$, such that 
for $0 \leq s \leq t \leq A$:
$$\mathbb{E}_{\mathbb{P}} \left[|X_t - X_s|^{\alpha} \right]  \leq C(A) (t-s)^{1 + \beta}.$$
Then $(X_t)_{t \geq 0}$ admits an adapted version which is locally H\"older of any index strictly smaller 
than $\beta / \alpha$. This version is unique up to indistinguishability.
\end{proposition}
\begin{proof}
Let $D$ be a countable, dense subset of $\mathbb{R}_+$. 
By the classical Kolmogorov criteria, there exists a version $Y$ of $X$ which is 
locally H\"older of any index strictly smaller 
than $\beta / \alpha$. 
One deduces that except on an N-negligible set, 
the restriction of $X$ to $D$ is also locally H\"older of any index strictly smaller 
than $\beta / \alpha$. One then extends this restriction to all $\mathbb{R}_+$ by continuity, 
which gives an adapted and H\"older version of $X$. 
\end{proof}
\noindent
One also has  a version of the d\'ebut theorem:
\begin{proposition} \label{debut}
Let $(\Omega, 
\mathcal{F}, (\mathcal{F}_t)_{t \geq 0}, \mathbb{P})$ be a filtered probability space satisfying the
N-usual conditions, 
and let $A$ be a progressive subset of $\mathbb{R}_+ \times \Omega$. Then the d\'ebut of $A$, i.e. 
the random time $D(A)$ such that for all $\omega \in \Omega$:
$$D(A)(\omega) := \inf \{t \geq 0, (t,\omega) \in A \}$$
is an $(\mathcal{F}_t)_{t \geq 0}$-stopping time. 
\end{proposition}
\begin{proof} 
Let $t \geq 0$. The set
$$A^{(t)} := A \, \cap \, ([0,t] \times \Omega)$$ 
is a progressive set of $[0,t] \times \Omega$, with respect to the filtered probability space
$(\Omega, \mathcal{F}_t, (\mathcal{F}_s)_{0 \leq s \leq t}, \mathbb{P})$.
Since this space satisfies the usual conditions, one can apply the classical d\'ebut theorem, which implies
that the d\'ebut $D(A^{(t)})$ of $A^{(t)}$ is a stopping time. Now, one immediately checks that 
$D(A) < t$ iff $D(A^{(t)}) < t$, which is in $\mathcal{F}_t$. Since $(\mathcal{F}_t)_{t \geq 0}$ is 
right-continuous, $D(A)$ is a stopping time. 
\end{proof}
\noindent
One has the following corollary:
\begin{corollary}
Let $(X_t)_{t \geq 0}$ be a progressively measurable process defined on a filtered probability 
space which satisfies the N-usual conditions 
(the condition of progressive mesurability is satisfied, in particular, if $(X_t)_{t \geq 0}$ is 
adapted and c\`adl\`ag). 
 Then, for all Borel sets $A \subset \mathbb{R}$:
$$T_A := \inf\{t \geq 0, X_t \in A\}$$
is an $(\mathcal{F}_t)_{t \geq 0}$-stopping time. 
\end{corollary}
\begin{proof}
The random time $T_A$ is the d\'ebut of the set:
$$\{(t, \omega) \in \mathbb{R}_+ \times \Omega, \, X_t(\omega) \in A\},$$
 which is progressive, since $(X_t)_{t \geq 0}$ is progressively measurable. 
\end{proof}
\noindent
Under the N-usual conditions, one can also prove the existence of the Doob-Meyer decomposition for submartingales.
In particular this implies the existence of a c\`adl\`ag (in fact continuous) and adapted version of the 
Brownian local time at level zero.  
\begin{proposition}
Let $(X_t)_{t \geq 0}$ be a right-continuous submartingale defined on a filtered probability
space $(\Omega, \mathcal{F}, (\mathcal{F}_t)_{t \geq 0}, \mathbb{P})$, satisfying the N-usual conditions.
 We suppose that $(X_t)_{t \geq 0}$ is of class $(DL)$, i.e. for all $a \geq 0$, $(X_T)_{T \in \mathcal{T}_a}$
is uniformly integrable, where $\mathcal{T}_a$ is the family of the $(\mathcal{F}_t)_{t \geq 0}$-stopping
times which are bounded by $a$ (for example, every nonnegative submartingale is of class $(DL)$).
Then, there exists a right-continuous $(\mathcal{F}_t)_{t \geq 0}$-martingale $(M_t)_{t \geq 0}$ and 
an increasing process
 $(A_t)_{t \geq 0}$ starting at zero, such that: $$X_t = M_t + A_t$$ for all $t \geq 0$, and for every
bounded, right-continuous martingale $(\xi_s)_{s \geq 0}$,
$$\mathbb{E} \left[\xi_t A_t \right] =\mathbb{E} \left[ \int_{(0,t]} \xi_{s-} dA_s \right],$$
where $\xi_{s-}$ is the left-limit of $\xi$ at $s$, almost surely well-defined for all $s > 0$.
The processes $(M_t)_{t \geq 0}$ and $(A_t)_{t \geq 0}$ are uniquely determined, up to indistinguishability. 
Moreover, they can be chosen to be continuous if $(X_t)_{t \geq 0}$ is a continuous process. 
\end{proposition}
\begin{proof}
Since for all $t \geq 0$, $(\Omega, \mathcal{F}_t, (\mathcal{F}_s)_{0 \leq s \leq t}, \mathbb{P})$
satisfies the usual conditions, there exists $(M^{(t)}_s)_{0 \leq s \leq t}$, a right-continuous
$(\mathcal{F}_s)_{0 \leq s \leq t}$-martingale and
 $(A^{(t)}_s)_{0 \leq s \leq t}$ increasing, such that $X_s = M^{(t)}_s + A^{(t)}_s$ for all $s \leq t$, 
and for all $u \in [0,t]$, and all bounded, right-continuous martingales $(\xi_s)_{s \geq 0}$:
$$\mathbb{E} \left[\xi_u A^{(t)}_u \right] =\mathbb{E} \left[ \int_{(0,u]} \xi_{s-} dA^{(t)}_s \right].$$
Moreover, one can suppose that $(M^{(t)}_s)_{0 \leq s \leq t}$ and $(A^{(t)}_s)_{0 \leq s \leq t}$ are
continuous if $(X_s)_{s \geq 0}$ is continuous. 
By uniqueness of $(M^{(t)}_s)_{0 \leq s \leq t}$ there exists an N-negligible set $N$ 
such that for all $\omega \notin N$, for all
integers $n \geq m \geq 1$, and for all $s \in [0, m]$: 
$$M^{(n)}_s (\omega) = M^{(m)}_s (\omega)$$
and
$$A^{(n)}_s (\omega) = A^{(m)}_s (\omega).$$
Since $N$ is N-negligible, $N \in \mathcal{F}_0$ with 
$\mathbb{P} [N]= 0$.
One can now define the processes $(M_s)_{s \geq 0}$ and $(A_s)_{s \geq 0}$ by
$M_s = M^{(n)}_s$, $A_s = A^{(n)}_s$ 
for $n \geq s \geq 0$, on the complement of $N$, and by $M_s = X_s$, $A_s = 0$ for all
 $s \geq 0$, on $N$. In particular, $(M_s)_{s \geq 0}$ and $(A_s)_{s \geq 0}$ are continuous
if $(X_s)_{s \geq 0}$ is continous.
\end{proof}
\noindent
One deduces the following corollary (giving quadratic variation):
\begin{corollary}
Let $(M_t)_{t \geq 0}$ be a continuous, square-integrable martingale defined on a filtered probability
space $(\Omega, \mathcal{F}, (\mathcal{F}_t)_{t \geq 0}, \mathbb{P})$, satisfying the N-usual conditions.
Then, there exists a unique continuous, increasing process $(\langle M \rangle_t)_{t \geq 0}$, starting from zero, 
 such that $(M_t^2 - \langle M \rangle_t)_{t \geq 0}$ is a martingale.
\end{corollary}
\noindent
The quadratic variation is involved in a very important way in the construction of the stochastic integral,
which can also be  made under N-usual conditions. We do not give here the details 
of the different constructions of the stochastic integral, but the general way to
go from the usual to the N-usual conditions is the following: let us suppose that under the usual conditions, one can 
construct a stochastic integral of the form:
$$\left( \int_0^s H_u dX_u \right)_{s\geq 0}$$
as a c\`adl\`ag, adapted process. Then, under the N-usual conditions, one can also define
$$\left( \int_0^s H_u dX_u \right)_{0 \leq s \leq t},$$
as a c\`adl\`ag, adapted process $(I^{(t)}_s)_{0 \leq s \leq t}$, since the restriction to $[0,t]$
of the underlying filtration satisfies the usual conditions. Now, let us assume that for $t'> t >0$, 
the restriction of $(I^{(t')}_s)_{0 \leq s \leq t'}$ to the interval $[0,t]$ is indistinguishable 
from the process $(I^{(t)}_s)_{0 \leq s \leq t}$. Since the processes are c\`adl\`ag, this assumption is 
implied by the fact that for all $s,t,t'$ such that $t' \geq t \geq s \geq 0$, 
$I^{(t)}_s = I^{(t')}_s$ almost surely, which is reasonable since $I^{(t)}_s$ and $I^{(t')}_s$ correspond 
to two constructions of the same stochastic integral:
$$\int_{0}^s H_u dX_u.$$
If this almost sure equality can be proved rigorously (this can be checked in each specific construction of the 
stochastic integral), one deduces that except on a N-negligible set $N$, all the processes 
$(I^{(n)}_s)_{0 \leq s \leq n}$, for $n \in \mathbb{N}$, are restrictions of each other.
By defining, for all $s \geq 0$, $I_s$ as zero on $N$, and as $I_s^{(n)}$ for $n \geq s$ on the complement of 
$N$, one obtains a c\`adl\`ag, adapted process, such that for $t \geq s \geq 0$:
$$I_s = I^{(t)}_s$$
almost surely, which gives a construction of the stochastic integral under the N-usual conditions.

Two other important results are the section theorem and the projection theorem. These results are 
proved in \cite{dell-meyer-projection}
 if the probability space is complete. 
In fact, they remain true if the space is only supposed to be N-complete. More precisely, the 
section theorem can be stated as follows:
\begin{proposition}
Let $(\Omega, \mathcal{F}, (\mathcal{F}_t)_{t \geq 0}, \mathbb{P})$ be an N-complete filtered probability 
space. Let $A$ be a optional subset of $\mathbb{R}_+ \times \Omega$. Then 
the image $\pi(A)$ of $A$ by the projection from $\mathbb{R}_+ \times \Omega$ to $\Omega$ is in $\mathcal{F}$,
and for all $\epsilon > 0$, there exists a stopping time $T$ enjoying the  following two properties:
\begin{itemize}
\item For all $\omega \in \Omega$ such that $T(\omega) < \infty$, $(T(\omega), \omega) \in A$.
\item $\mathbb{P} [T < \infty] \geq \mathbb{P} [\pi(A)] - \epsilon$.
\end{itemize}
\end{proposition} 
\begin{proof}
Let $\epsilon > 0$ and $u \in \mathbb{R}_+$. The negligible sets of $(\Omega, \mathcal{F}_u, \mathbb{P})$ are 
N-negligible with respect to $(\Omega, \mathcal{F}, (\mathcal{F}_t)_{t \geq 0}, \mathbb{P})$,
hence, they are in $\mathcal{F}_0$: the probability space
 $(\Omega, \mathcal{F}_u, \mathbb{P})$ is complete. 
Moreover, the set $A \cap \left([0,u) \times \Omega \right)$ is optional with respect to 
$(\Omega, \mathcal{F}_u, (\mathcal{F}_{t \wedge u})_{t \geq 0}, \mathbb{P})$. This is a consequence of 
the  following two facts:
\begin{itemize}
\item The family of sets $B \subset \mathbb{R}_+ \times \Omega$ such that $B \cap \left([0,u) \times \Omega
 \right)$ 
is optional 
with respect to $(\Omega, \mathcal{F}_u, (\mathcal{F}_{t \wedge u})_{t \geq 0}, \mathbb{P})$ is a 
$\sigma$-algebra;
\item For all stopping times $T$ with respect to $(\Omega, \mathcal{F}, (\mathcal{F}_t)_{t \geq 0}, \mathbb{P})$,
the set of $(t, \omega) \in [0,u) \times \Omega$ such that $T(\omega) \leq t$ is optional 
with respect to $(\Omega, \mathcal{F}_u, (\mathcal{F}_{t \wedge u})_{t \geq 0}, \mathbb{P})$.
\end{itemize}
\noindent
The first property is proved as follows:
\begin{itemize}
\item If $B_n \cap \left( [0,u) \times \Omega \right)$ is optional for all integers $n \geq 1$:
$$ \left(\bigcup_{n \geq 1} B_n \right) \, \cap \, \left( [0,u) \times \Omega \right)
 = \bigcup_{n \geq 1} \, \left[B_n \cap  \left([0,u) \times \Omega \right) \, \right]$$
is optional. 
\item If $B \cap \left( [0,u) \times \Omega \right)$ is optional, and if $C$ is the complement of $B$, then 
$C \cap \left( [0,u) \times \Omega \right)$ is the intersection of $[0,u) \times \Omega$ and 
the complement of $B \cap \left( [0,u) \times \Omega \right)$, hence, it is optional. 
\end{itemize}
\noindent
The set involved in the second property corresponds to the stochastic interval $[T,u)$ for 
$T < u$, and the empty set for $T \geq u$. Equivalently it corresponds to the intersection of 
$[T \wedge u, \infty)$ and $[0,u)$, and then is optional since $T \wedge u$ is a stopping time with respect
to the filtration $(\mathcal{F}_{t \wedge u})_{t \geq 0}$. 
We have now proved that for all integers $n \geq 1$,
 $A_n := A \cap \left([0,n) \times \Omega \right)$ is optional with respect to 
$(\Omega, \mathcal{F}_n, (\mathcal{F}_{t \wedge n})_{t \geq 0}, \mathbb{P})$. By the section theorem in 
\cite{dell-meyer-projection}, $A_n$ is measurable and there exists an $(\mathcal{F}_{t \wedge n})_{t \geq 0}$-stopping time
$T_n$ such that $T_n ( \omega ) < \infty$ implies $(T_n(\omega), \omega) \in A_n$ and 
$\mathbb{P} [T_n < \infty] \geq \mathbb{P} [\pi(A_n)] - (\epsilon/2^n)$, or equivalently, 
$$\mathbb{P} [\pi(A_n) \backslash \, \{T_n < \infty \} \, ] \leq \epsilon/2^n,$$
since $\pi(A_n)$ contains the event $\{T_n < \infty\}$.  
Now $\pi(A)$ is the increasing union of $\pi(A_n)$ for $n \geq 1$, hence, it is in $\mathcal{F}$. 
Moreover there exists $N \geq 1$ such that $\mathbb{P} [\pi(A_N)] \geq \mathbb{P} [\pi(A)] - \epsilon$. 
Let $T$ be the infimum of $T_n$ for $1 \leq n \leq N$: it is an $(\mathcal{F}_t)_{t \geq 0}$-stopping 
time. If $T$ is finite, let $n$ be an index such that $T = T_n$; one 
has $(T_n(\omega), \omega) \in A_n$, which implies that $(T(\omega), \omega) \in A$. 
Moreover
\begin{align*}
\mathbb{P} [\pi(A_N) \, \backslash \, \{T < \infty\} ] & = 
\mathbb{P} \left[ \left( \bigcup_{1 \leq n \leq N} \, \pi(A_n) \right) \, \cap \, 
\left( \bigcap_{1 \leq n \leq N} \, \{T_n = \infty\} \, \right) \right]
\\ & \leq \sum_{1 \leq n \leq N} \mathbb{P} [\pi(A_n) \cap \{T_n = \infty\} ] \\ 
 & \leq \epsilon.
\end{align*}
\noindent
Therefore
$$\mathbb{P} [T < \infty] \geq \mathbb{P}[\pi(A_N)] - \epsilon \geq \mathbb{P}[A] - 2 \epsilon.$$
\end{proof}
\noindent
We deduce from the section theorem the following version of the projection theorem:
\begin{proposition}
Let $(X_t)_{t \geq 0}$ be a bounded, measurable process, defined on a N-complete, filtered probability 
space $(\Omega, \mathcal{F}, (\mathcal{F}_t)_{t \geq 0}, \mathbb{P})$. Then, there exists an optional 
stochastic process $(Y_t)_{t \geq 0}$ such that for all stopping times $T$:
$$\mathbb{E} [X_T \mathds{1}_{T < \infty} \, | \mathcal{F}_T ] = Y_T \mathds{1}_{T < \infty}.$$
This process is unique up to indistinguishability. 
\end{proposition}
\begin{proof}
The proof in \cite{dell-meyer-projection} is sufficient here, since it uses completeness of the space only via the section theorem. 
\end{proof}

\section{Extension of measures and properties preserved by the N-augmentation} \label{extensions}
\noindent
A natural question which arises when one has, on a filtered probability space $(\Omega, \mathcal{F}, 
(\mathcal{F}_t)_{t \geq 0}, \mathbb{P})$, a coherent family $(\mathbb{Q}_t)_{t \geq 0}$ of
probability measures, $\mathbb{Q}_t$ defined on $\mathcal{F}_t$, is the following:
does there exist a probability measure $\mathbb{Q}$ defined on $\mathcal{F}$, such that its restriction to 
$\mathcal{F}_t$ is $\mathbb{Q}_t$ for all $t \geq 0$? As we have seen in the introduction, 
the answer to this question is not obvious at all and can be negative in very simple cases.
However, one can state sufficient conditions on the space
$(\Omega, \mathcal{F}, 
(\mathcal{F}_t)_{t \geq 0}, \mathbb{P})$, under which the answer is always positive. This problem is not new and has already received attention in the literature: it is dealt with in great details by Parthasarathy in \cite{Parth}. In particular the conditions we are going to state are not new and are already contained in \cite{Parth}, p. 141. We nevertheless provide a proof, first because our method is slightly different and uses the existence of a regular conditional probability measure to construct the measure on $\mathcal{F}_\infty$, and second because this type of result does not seem to be well known. To be complete, one should add that H. F\"ollmer, in \cite{Foellmer},  has also studied a similar problem in a more general context  and has obtained sufficient conditions relying on Parthasarathy's results. However the results there would not  apply here: in particular the conditions  in \cite{Foellmer}  are not satisfied by the space of continuous functions and  the space of  c\`adl\`ag functions, equipped with the natural filtration of the canonical process (these filtrations are not standard filtrations with the definition of a standard filtration given in \cite{Foellmer}).

We first start with some definitions and then give sufficient conditions on a filtered measurable space $(\Omega, \mathcal{F}, (\mathcal{F}_t)_{t \geq 0})$ in order to have  a positive answer to the extension of measure problem. Then we show that when this extension is possible, then it still holds under the N-usual augmentation with respect to some probability measure. Whenever we need a reference for probability measures on metric spaces, we cite \cite{Parth} or the first chapter of \cite{SV}.

\begin{definition} \label{P}
Let $(\Omega, \mathcal{F}, (\mathcal{F}_t)_{t \geq 0})$ be a filtered measurable space, such that
$\mathcal{F}$ is the $\sigma$-algebra generated by $\mathcal{F}_t$, 
$t \geq 0$: $\mathcal{F}=\bigvee_{t\geq0}\mathcal{F}_t$. We shall say that the
 property (P) holds if and only if $(\mathcal{F}_t)_{t \geq 0}$ 
enjoys the following conditions: 
\begin{itemize}
\item For all $t \geq 0$, $\mathcal{F}_t$ is generated by a countable number of sets;
\item For all $t \geq 0$, there exists a Polish space $\Omega_t$, and a surjective map 
 $\pi_t$ from $\Omega$ to $\Omega_t$, such that $\mathcal{F}_t$ is the $\sigma$-algebra of the inverse
 images, by $\pi_t$, of Borel sets in $\Omega_t$, and such that for all $B \in \mathcal{F}_t$, 
 $\omega \in \Omega$, $\pi_t (\omega) \in \pi_t(B)$ implies $\omega \in B$;
\item If $(\omega_n)_{n \geq 0}$ is a sequence of elements of $\Omega$, such that for all $N \geq 0$,
$$\bigcap_{n = 0}^{N} A_n (\omega_n) \neq \emptyset,$$
where $A_n (\omega_n)$ is the intersection of the sets in $\mathcal{F}_n$ containing $\omega_n$, 
then:
$$\bigcap_{n = 0}^{\infty} A_n (\omega_n) \neq \emptyset.$$
\end{itemize}
\end{definition}
\noindent
Given this technical definition, one can state the following result:
\begin{proposition} \label{extension}
Let $(\Omega, \mathcal{F}, (\mathcal{F}_t)_{t \geq 0})$ be a filtered measurable space satisfying the 
property (P), and let, for $t \geq 0$, $\mathbb{Q}_t$ be a probability measure on $(\Omega, \mathcal{F}_t)$, 
such that for all $t \geq s \geq 0$, $\mathbb{Q}_s$ is the restriction of $\mathbb{Q}_t$ to $\mathcal{F}_s$.
Then, there exists a unique measure $\mathbb{Q}$ on $(\Omega, \mathcal{F})$ such that for all $t \geq 0$,
its restriction to $\mathcal{F}_t$ is equal to $\mathbb{Q}_t$. 
\end{proposition}
\begin{proof} The uniqueness is a direct application of the monotone class theorem. Let us prove the 
existence. 
For $t \geq 0$, the map $\pi_t$ is by assumption measurable from $(\Omega, \mathcal{F}_t)$ 
to $(\Omega_t, \mathcal{B} (\Omega_t))$ (where $\mathcal{B} (\Omega_t)$ is the Borel $\sigma$-algebra
of $\Omega_t$). Let $\bar{\mathbb{Q}}_t$ be the image of $\mathbb{Q}_t$ by 
$\pi_t$. By Theorem 1.1.6 of \cite{SV}, for $0 \leq s \leq t$, there exists a conditional probability distribution 
of $\bar{\mathbb{Q}}_t$ given the $\sigma$-algebra $\pi_t(\mathcal{F}_s)$, generated by the 
images by $\pi_t$ of the sets in $\mathcal{F}_s$. Note that this $\sigma$-algebra
 is included in $\mathcal{B} (\Omega_t)$. Indeed, if $B \in \mathcal{F}_s$, there exists 
$A \in \mathcal{B} (\Omega_t)$ such that $B = \pi_t^{-1} (A)$, and then 
$\pi_t (B) = \pi_t \circ \pi_t^{-1}(A)$, which is equal to $A$ by the surjectivity of $\pi_t$. 
 Now, the existence of the conditional probability distribution described above means that one can find
 a family $(Q_{\omega})_{\omega \in \Omega_t}$ of probability measures on $(\Omega_t, \mathcal{B}
(\Omega_t))$ such that:
\begin{itemize}
\item For each $B \in \mathcal{B}(\Omega_t)$, $\omega \rightarrow Q_{\omega} (B)$ is 
$\pi_t (\mathcal{F}_s)$-measurable;
\item For every $A \in \pi_t (\mathcal{F}_s)$, $B \in \mathcal{B} (\Omega_t)$:
$$\bar{\mathbb{Q}}_t (A \cap B) = \int_A Q_{\omega} (B) \, \bar{\mathbb{Q}}_t(d \omega).$$ 
\end{itemize}
\noindent
Now, for all $\omega \in \Omega$, let us define the map $R_{\omega}$ from $\mathcal{F}_t$ to 
$\mathbb{R}_+$, by
$$R_{\omega} [B] =  Q_{\pi_t (\omega)} [ \pi_t (B)].$$
The map $R_{\omega}$ is a probability measure on $(\Omega, \mathcal{F}_t)$. Indeed, 
$$R_{\omega} [\Omega] = Q_{\pi_t (\omega)}  [\Omega_t] = 1,$$
since $\pi_t$ is surjective. Moreover, let $(B_k)_{k \geq 1}$ be a family of disjoint sets in 
$\mathcal{F}_t$, and $B_0$ their union. By assumption, there exist $(\tilde{B}_k)_{k \geq 0}$ in $\mathcal{B}
 (\Omega_t)$
such that $B_k = \pi_t^{-1} (\tilde{B}_k)$. Since $\pi_t$ is surjective, $\pi_t (B_k) = \tilde{B}_k$. 
Moreover, the sets $(\tilde{B}_k)_{k \geq 1}$ are pairwise disjoint. Indeed, if $x \in  
\tilde{B}_k \cap \tilde{B}_l$ for $k > l \geq 1$, then, by surjectivity, there exists $y \in \Omega$,
such that $x = \pi_t(y)$, which implies $y \in \pi_t^{-1} (\tilde{B}_k) \cap \pi_t^{-1} (\tilde{B}_l)$,
and then $y \in B_k \cap B_l$, which is impossible. Therefore:
\begin{align*}
R_{\omega} [B_0] & =  Q_{\pi_t (\omega)} [ \tilde{B}_0] \\ & = 
\sum_{k \geq 1} Q_{\pi_t (\omega)} [ \tilde{B}_k] \\ & = 
\sum_{k \geq 1} R_{\omega} [B_k]. 
\end{align*} 
\noindent
Hence, $R_{\omega}$ is a probability measure. Moreover,
 for each $B \in \mathcal{F}_t$, the map $\omega \rightarrow R_{\omega} (B)$ is the composition
of the measurable maps $\omega \rightarrow  \pi_t(\omega)$ from $(\Omega, \mathcal{F}_s)$
 to $(\Omega_t, \pi_t (\mathcal{F}_s))$, and $\omega' \rightarrow 
Q_{\omega'}[\pi_t(B)]$ from $(\Omega_t, \pi_t (\mathcal{F}_s))$ to $(\mathbb{R}_+, \mathcal{B} (\mathbb{R}_+))$,
and hence, it is $\mathcal{F}_s$-measurable. The measurability of $\pi_t$ follows from the fact that 
for $A \in \mathcal{F}_s$, the inverse image of $\pi_t(A)$ by $\pi_t$ is exactly $A$ (by an assumption
given in the definition of the property (P)). Moreover, for every $A \in  \mathcal{F}_s$, $B \in
\mathcal{F}_t$:
\begin{align*}
\mathbb{Q}_t (A \cap B) & = \mathbb{Q}_t [(\pi_t^{-1} \circ \pi_t (A)) 
 \cap (\pi_t^{-1} \circ \pi_t (B))] \\ & =  \mathbb{Q}_t  [\pi_t^{-1} 
( \pi_t(A) \cap \pi_t(B))] \\ & = \bar{\mathbb{Q}}_t [\pi_t(A) \cap \pi_t(B)]
\\ & = \int_{\Omega_t} \mathds{1}_{\omega \in \pi_t(A)}
Q_{\omega} (\pi_t (B)) \bar{\mathbb{Q}}_t(  d \omega) \\
& = \int_{\Omega} \mathds{1}_{\pi_t(\omega) \in \pi_t(A)} \, Q_{\pi_t(\omega)} (\pi_t (B)) 
\mathbb{Q}_t (d \omega) \\ & = \int_A R_{\omega} (B) \mathbb{Q}_t (d \omega)
\end{align*}
\noindent
Finally, we have found a conditional probability distribution of $\mathbb{Q}_t$
with respect to $\mathcal{F}_s$. Since $\mathcal{F}_s$ is countably generated, this conditional 
probability distribution is regular, by Theorem 1.1.8 in \cite{SV}. One can then apply Theorem
1.1.9, again in \cite{SV}, and since $\mathcal{F}$ is the 
$\sigma$-algebra generated by $\mathcal{F}_t$, $t \geq 0$, 
one obtains a probability distribution $\mathbb{Q}$ on $(\Omega, \mathcal{F})$ such that 
for all integers $n \geq 0$, the restriction of $\mathbb{Q}$ to $\mathcal{F}_n$ is 
$\mathbb{Q}_n$. Now, for $t \geq 0$, let $\Lambda_t$ be an event in $\mathcal{F}_t$. 
One has, for $n > t$, integer:
$$\mathbb{Q} [\Lambda_t] = \mathbb{Q}_n [\Lambda_t] =  \mathbb{Q}_t[\Lambda_t],$$
which implies that $\mathbb{Q}$ satisfies the assumptions of Proposition \ref{extension}.
\end{proof}
\noindent
Now, it remains to find a way to see if a space satisfies property (P) or not. We do not give a general 
condition here, but we prove that it is the case for the most classical spaces, 
endowed with their canonial filtration.
More precisely:
\begin{proposition}
Let 	 $\Omega$ be $\mathcal{C}(\mathbb{R}_+,\mathbb{R}^d)$, the space of continuous functions from 
	$\mathbb{R}_+$ to $\mathbb{R}^d$, or $\mathcal{D}(\mathbb{R}_+,\mathbb{R}^d)$, the space of c\`adl\`ag functions
 from $\mathbb{R}_+$ 
	to $\mathbb{R}^d$ (for some $d \geq 1$). For $t \geq 0$, define $(\mathcal{F}_t)_{t \geq 0}$
as the natural filtration of the canonical process $Y$, 
and $\mathcal{F}=\bigvee_{t\geq0}\mathcal{F}_t$. Then  $(\Omega, \mathcal{F},(\mathcal{F}_t)_{t \geq 0})$
 satisfies property (P).
\end{proposition}
\begin{proof}
  Let us prove
this result for c\`adl\`ag functions (for continuous functions, the result is similar and
  is proved in \cite{SV}). For all $t \geq 0$, $\mathcal{F}_t$ is 
generated by the variables $Y_{rt}$, for $r$, rational, in $[0,1]$, hence, it is countably generated.
For the second property, one can take for $\Omega_t$, the set of c\`adl\`ag functions
from $[0,t]$ to $\mathbb{R}^d$, and for $\pi_t$, the restriction to the interval $[0,t]$.
The space $\Omega_t$ is Polish if one endows it with the Skorokhod metric. Moreover, 
its Borel $\sigma$-algebra is equal to the $\sigma$-algebra generated by the coordinates, a result 
from which one easily  deduces the properties of $\pi_t$ which need to be satisfied. The 
third property is easy to check: let us suppose that $(\omega_n)_{n \geq 0}$ is a sequence
 of elements of $\Omega$, such that for all $N \geq 0$,
$$\bigcap_{n = 0}^{N} A_n (\omega_n) \neq \emptyset,$$
where $A_n (\omega_n)$ is the intersection of the sets in $\mathcal{F}_n$ containing $\omega_n$.
Here,  $A_n (\omega_n)$ is the set of functions $\omega'$ which coincide with 
$\omega_n$ on $[0,n]$. Moreover, for $n \leq n'$, integers, the intersection
of $A_n (\omega_n)$ and $A_{n'} (\omega_{n'})$ is not empty, and then $\omega_n$ and $\omega_{n'}$ coincide 
on $[0,n]$. Therefore, there exists a c\`adl\`ag function $\omega$ which coincides
with $\omega_n$ on $[0,n]$, for all $n$, which implies:
$$\bigcap_{n = 0}^{\infty} A_n (\omega_n) \neq \emptyset.$$
\end{proof}
\begin{remark}
It is easily seen that the conditions of Proposition \ref{extension} are not satisfied by the space $\mathcal{C}([0,1],\mathbb{R})$ endowed with the filtration generated by the canonical process; an explicit counter example is provided in \cite{FI}.
\end{remark}
 \noindent
We have now some examples of filtered probability spaces for which the extension of measures 
are always possible. However, these spaces do not possess any property of completion, and their filtration
is not right-continuous. That is why we need to check that, at least in many interesting cases, 
this theorem of extension is still available after taking the N-augmentation.
Before giving the corresponding statement, let us give a simple result, proving 
that many properties of stochastic 
processes are not changed by the N-augmentation of the underlying filtration.
\begin{proposition} \label{conditional}
Let  $(\Omega, \mathcal{F}, (\mathcal{F}_t)_{t \geq 0}, \mathbb{P})$ be a filtered probability space, 
and $(\Omega, \widetilde{\mathcal{F}}, (\widetilde{\mathcal{F}}_t)_{t \geq 0}, \widetilde{\mathbb{P}})$
its N-augmentation. Let $X$ be an $\mathcal{F}$-measurable random variable, integrable with respect
to $\mathbb{P}$. Then $X$ is also integrable with respect to $\widetilde{\mathbb{P}}$, and 
$$\mathbb{E}_{\widetilde{\mathbb{P}}} [X] = \mathbb{E}_{\mathbb{P}}[X].$$
Moreover, for all $t \geq 0$
$$\mathbb{E}_{\widetilde{\mathbb{P}}} [X| \widetilde{\mathcal{F}}_t] =
\mathbb{E}_{\mathbb{P}} [X|\mathcal{F}_t],$$
$\widetilde{\mathbb{P}}$-almost surely.
\end{proposition}
\begin{proof}
For all sets $\Lambda \in \mathcal{F}$, one has
$$\widetilde{\mathbb{P}} [\Lambda] = \mathbb{P}[\Lambda].$$
One deduces that 
$$\mathbb{E}_{\widetilde{\mathbb{P}}} [X] = \mathbb{E}_{\mathbb{P}}[X],$$
for every nonnegative, $\mathcal{F}$-measurable random variable $X$, and then for every integrable variable 
$X$. Let $t \geq 0$ and let $\Lambda_t$ be an event in $\widetilde{\mathcal{F}}_t$. There exists
 an event $\Lambda'_t
\in \mathcal{F}_t$ such that its symmetric difference with $\Lambda_t$ is N-negligible.
One then has
\begin{align*}
\mathbb{E}_{\widetilde{\mathbb{P}}}[X \, \mathds{1}_{\Lambda_t} ]
= \mathbb{E}_{\widetilde{\mathbb{P}}} [X \, \mathds{1}_{\Lambda'_t}]
 & = \mathbb{E}_{\mathbb{P}} [X \, \mathds{1}_{\Lambda'_t}]
 = \mathbb{E}_{\mathbb{P}} [ \mathbb{E}_{\mathbb{P}}  [X |\mathcal{F}_t] \, \mathds{1}_{\Lambda'_t} ]
\\ & = \mathbb{E}_{\widetilde{\mathbb{P}}} [ \mathbb{E}_{\mathbb{P}} 
 [X |\mathcal{F}_t] \, \mathds{1}_{\Lambda'_t} ]
 = \mathbb{E}_{\widetilde{\mathbb{P}}} [ \mathbb{E}_{\mathbb{P}} 
 [X |\mathcal{F}_t] \, \mathds{1}_{\Lambda_t} ],
\end{align*}
\noindent
which proves Proposition \ref{conditional}. 
\end{proof}
\noindent
Proposition \ref{conditional} has, in particular, the following consequence:
\begin{corollary} \label{conditional2}
Let  $(\Omega, \mathcal{F}, (\mathcal{F}_t)_{t \geq 0}, \mathbb{P})$ be a filtered probability space, 
and $(\Omega, \widetilde{\mathcal{F}}, (\widetilde{\mathcal{F}}_t)_{t \geq 0}, \widetilde{\mathbb{P}})$
its N-augmentation. Then if $(X_t)_{t \geq 0}$ is a (sub)martingale with respect
to the filtration $(\mathcal{F}_t)_{t \geq 0}$ and the probability $\mathbb{P}$, then 
it is also a (sub)martingale with respect to  $(\widetilde{\mathcal{F}}_t)_{t \geq 0}$  and 
$\widetilde{\mathbb{P}}$.
\end{corollary}
\begin{proof}
Let us suppose that $(X_t)_{t \geq 0}$ is a submartingale with respect 
to $(\mathcal{F}_t)_{t \geq 0}$ and $\mathbb{P}$.
By Proposition \ref{conditional}, for all $t \geq 0$, $X_t$ is integrable with respect 
to $\widetilde{\mathbb{P}}$, and for $s \leq t$:
$$\mathbb{E}_{\widetilde{\mathbb{P}}} [X_t |\widetilde{\mathcal{F}}_s]=
\mathbb{E}_{\mathbb{P}} [X_t| \mathcal{F}_s] \geq X_s,$$
which implies that $(X_t)_{t \geq 0}$ is a submartingale with respect to 
$(\widetilde{\mathcal{F}}_t)_{t \geq 0}$  and  $\widetilde{\mathbb{P}}$. The proof 
for martingales is exactly similar. 
\end{proof}
\noindent
From Corollary \ref{conditional2}, one deduces:
\begin{corollary}
Let  $(\Omega, \mathcal{F}, (\mathcal{F}_t)_{t \geq 0}, \mathbb{P})$ be a filtered probability space, 
and $(\Omega, \widetilde{\mathcal{F}}, (\widetilde{\mathcal{F}}_t)_{t \geq 0}, \widetilde{\mathbb{P}})$
its N-augmentation. Then if $(X_t)_{t \geq 0}$ is a standard Brownian motion with respect 
to the filtration $(\mathcal{F}_t)_{t \geq 0}$ and the probability $\mathbb{P}$, then 
it is also a Brownian motion with respect to  $(\widetilde{\mathcal{F}}_t)_{t \geq 0}$  and 
$\widetilde{\mathbb{P}}$.
\end{corollary}
\begin{proof}
We apply Corollary \ref{conditional2} to $(X_t)_{t \geq 0}$ and $(X_t^2 - t)_{t \geq 0}$, and by continuity
of $(X_t)_{t \geq 0}$, we are done. 
\end{proof}
\noindent
Note that the results above are still available under usual augmentation, instead 
of N-augmentation (recall that usual 
augmentation is obtained by putting all the $\mathbb{P}$-negligible sets in 
$\mathcal{F}$ and $\mathcal{F}_t$ for 
all $t \geq 0$, and by completing the probability $\mathbb{P}$). However, in order to conserve the
extension of coherent families of probabilties, on really needs to consider the N-augmentation instead
 of the usual augmentation. The precise statement is the following: 
\begin{proposition} \label{extensionaugmentation}
Let $(\Omega, \mathcal{F}, (\mathcal{F}_t)_{t \geq 0}, \mathbb{P})$ be a filtered probability space, 
and $(\Omega, \widetilde{\mathcal{F}}, (\widetilde{\mathcal{F}}_t)_{t \geq 0}, \widetilde{\mathbb{P}})$
its N-augmentation. We assume that for all coherent families of probability measures
$(\mathbb{Q}_t)_{t \geq 0}$, 
such that $\mathbb{Q}_t$ is defined on $\mathcal{F}_t$ for all $t \geq 0$,  
 there exists a unique probability measure $\mathbb{Q}$ on $\mathcal{F}$ which coincides with $\mathbb{Q}_t$
on $\mathcal{F}_t$ for all $t \geq 0$. Then, for all coherent families of probability measures
$(\widetilde{\mathbb{Q}}_t)_{t \geq 0}$, such that $\widetilde{\mathbb{Q}}_t$ is defined 
on $\widetilde{\mathcal{F}}_t$,
and is absolutely continuous with respect to the restriction
of $\widetilde{\mathbb{P}}$ to $\widetilde{\mathcal{F}}_t$, there exists a unique probability measure
$\widetilde{\mathbb{Q}}$
on $\widetilde{\mathcal{F}}$ which coincides with $\widetilde{\mathbb{Q}}_t$ on $\widetilde{\mathcal{F}}_t$ 
for all $t \geq 0$. 
\end{proposition}
\begin{proof}
Let $(\widetilde{\mathbb{Q}}_t)_{t \geq 0}$ be a coherent family of probability measures such that
$\widetilde{\mathbb{Q}}_t$ is defined on $\widetilde{\mathcal{F}}_t$ and is absolutely continuous with 
respect to the restriction of $\widetilde{\mathbb{P}}$ to $\widetilde{\mathcal{F}}_t$. For all 
$t \geq 0$, one can consider the restriction $\mathbb{Q}_t$ of $\widetilde{\mathbb{Q}}_t$ 
to $\mathcal{F}_t$, and the family $(\mathbb{Q}_t)_{t \geq 0}$ is coherent. Indeed, for all $s \in [0,t]$ and for 
all events $\Lambda_s \in \mathcal{F}_s$, one has:
$$\mathbb{Q}_t [\Lambda_s ] = \widetilde{\mathbb{Q}}_t [\Lambda_s]
= \widetilde{\mathbb{Q}}_s [\Lambda_s] = \mathbb{Q}_s [\Lambda_s].$$
Therefore there exists a unique probability measure $\mathbb{Q}$ on $\mathcal{F}$ such that for all $t \geq 0$, 
the restriction of $\mathbb{Q}$ to $\mathcal{F}_t$ is $\mathbb{Q}_t$. 
Now, let $(\Omega, \mathcal{F}', (\mathcal{F}'_t)_{t \geq 0}, \mathbb{Q}')$ be the N-augmentation
of $(\Omega, \mathcal{F}, (\mathcal{F}_t)_{t \geq 0}, \mathbb{Q})$. If an event $A$ is 
N-negligible with respect to $(\Omega, \mathcal{F}, (\mathcal{F}_t)_{t \geq 0}, \mathbb{P})$,
there exists $(B_n)_{n \geq 0}$, such that
$$A \subset \bigcup_{n \geq 0} B_n,$$
and for all $n \geq 0$, $B_n \in \mathcal{F}_n$ and $\mathbb{P} [B_n] = 0$. Since by assumption,
$\widetilde{\mathbb{Q}}_n$ is absolutely continuous with respect to the restriction of $\widetilde{\mathbb{P}}$
to $\widetilde{\mathcal{F}}_n$:
$$\mathbb{Q} [B_n] = \mathbb{Q}_n [B_n] = \widetilde{\mathbb{Q}}_n [B_n] = 
\widetilde{\mathbb{P}} [B_n] = \mathbb{P}[B_n]=0.$$
One deduces that $A$ is N-negligible with respect to 
$(\Omega, \mathcal{F}, (\mathcal{F}_t)_{t \geq 0}, \mathbb{Q})$. In other words, all the 
N-negligible sets of $(\Omega, \mathcal{F}, (\mathcal{F}_t)_{t \geq 0}, \mathbb{P})$ 
are also N-negligible with respect to $(\Omega, \mathcal{F}, (\mathcal{F}_t)_{t \geq 0}, \mathbb{Q})$,
which implies that $\mathcal{F}'$ contains $\widetilde{\mathcal{F}}$. 
Therefore, one can define the probability measure $\widetilde{\mathbb{Q}}$ as the 
restriction of $\mathbb{Q}'$ to the $\sigma$-algebra $\widetilde{\mathcal{F}}$: it remains to check
that this measure satisfies Proposition \ref{extensionaugmentation}. Indeed, let $\Lambda_t$ be 
an event in $\widetilde{\mathcal{F}}_t$, for some $t \geq 0$. There exists $\Lambda'_t 
\in \mathcal{F}_{t+}$ such that its symmetric difference with $\Lambda_t$ is N-negligible with respect
to $(\Omega, \mathcal{F}, (\mathcal{F}_t)_{t \geq 0}, \mathbb{P})$,
and then also with respect to $(\Omega, \mathcal{F}, (\mathcal{F}_t)_{t \geq 0}, \mathbb{Q})$.
One deduces that
\begin{align*}
\widetilde{\mathbb{Q}} [\Lambda_t] = \mathbb{Q}'[\Lambda_t] & = \mathbb{Q}[\Lambda'_t]
 = \mathbb{Q}_{t+1} [\Lambda'_t] \\ 
& = \widetilde{\mathbb{Q}}_{t+1} [\Lambda'_t] = \widetilde{\mathbb{Q}}_{t+1} [\Lambda_t] \\
& = \widetilde{\mathbb{Q}}_t [\Lambda_t]. 
\end{align*}
\noindent
Here the equality
$$\widetilde{\mathbb{Q}}_{t+1} [\Lambda'_t] = \widetilde{\mathbb{Q}}_{t+1} [\Lambda_t]$$
is due to the fact that $\Lambda'_t$ and $\Lambda_t$ are both in $\widetilde{\mathcal{F}}_{t+1}$, 
their symmetric difference is N-negligible with respect to
$(\Omega, \mathcal{F}, (\mathcal{F}_t)_{t \geq 0}, \mathbb{P})$, and $\widetilde{\mathbb{Q}}_{t+1}$ 
 is absolutely continuous with respect to the restriction of $\widetilde{\mathbb{P}}$ 
to $\widetilde{\mathcal{F}}_{t+1}$. The uniqueness of $\widetilde{\mathbb{Q}}$ is proved as follows:
its restriction $\mathbb{Q}$ to $\mathcal{F}$ is uniquely determined, since it has to coincide 
with $\mathbb{Q}_t$ on $\mathcal{F}_t$ (recall that $\mathbb{Q}_t$ is 
the restriction of $\widetilde{\mathbb{Q}}_t$ to $\mathcal{F}_t$). Now, for all events
$\Lambda \in \widetilde{\mathcal{F}}$, there exists $\Lambda' \in \mathcal{F}$ such that the symmetric 
difference of $\Lambda$ and $\Lambda'$ is N-negligible with respect to 
$(\Omega, \mathcal{F}, (\mathcal{F}_t)_{t \geq 0}, \mathbb{P})$, and then 
with respect to $(\Omega, \mathcal{F}, (\mathcal{F}_t)_{t \geq 0}, \mathbb{Q})$. One deduces 
that $(\Lambda \backslash \Lambda') \cup (\Lambda' \backslash \Lambda)$ is included
in a set $B \in \mathcal{F}$ such that $\mathbb{Q} [B] = 0$, and then
$$\widetilde{\mathbb{Q}} \left[\Lambda \backslash \Lambda') \cup (\Lambda' \backslash \Lambda) \right]
\leq \widetilde{\mathbb{Q}} [B] = \mathbb{Q} [B] = 0,$$
which implies
$$\widetilde{\mathbb{Q}} [\Lambda] = \widetilde{\mathbb{Q}} [\Lambda'] 
= \mathbb{Q} [\Lambda'],$$
and $\widetilde{\mathbb{Q}} [\Lambda]$ is uniquely determined.
\end{proof}
\noindent 
One immediately deduces, from Propositions \ref{extension} and \ref{extensionaugmentation}, the following:
\begin{corollary}
Let $(\Omega, \mathcal{F}, (\mathcal{F}_t)_{t \geq 0}, \mathbb{P})$ be the N-augmentation of 
a filtered probability space satisfying the property (P). Then if $(\mathbb{Q}_t)_{t \geq 0}$ 
is a coherent family of probability measures, $\mathbb{Q}_t$ defined on $\mathcal{F}_t$, and absolutely 
continuous with respect to the restriction of $\mathbb{P}$ to $\mathcal{F}_t$, there exists
a unique probability measure $\mathbb{Q}$ on $\mathcal{F}$ which coincides with $\mathbb{Q}_t$ on 
$\mathcal{F}_t$, for all $t \geq 0$. 
\end{corollary}
\noindent
We observe that if $(\Omega, \mathcal{F}, (\mathcal{F}_t)_{t \geq 0})$ satisfies the property (P), 
then $\mathcal{F}$ is the $\sigma$-algebra generated by $\mathcal{F}_t$, $t \geq 0$. This property is 
clearly preserved when one takes the N-augmentation.

\end{document}